\documentclass[10pt]{amsart}
\usepackage{latexsym, amsmath,amssymb}

\renewcommand{\epsilon}{\varepsilon}

\newtheorem{theorem}{Theorem}
\newtheorem{lemma}[theorem]{Lemma}
\newtheorem{corr}[theorem]{Corollary}

\newtheorem{proposition}[theorem]{Proposition}
\newtheorem{deff}[theorem]{Definition}

\newcommand{\bth}{\begin{theorem}}
\newcommand{\ble}{\begin{lemma}}
\newcommand{\bcor}{\begin{corr}}

\newcommand{\bdeff}{\begin{deff}}

\newcommand{\bprop}{\begin{proposition}}
\newcommand{\ele}{\end{lemma}}
\newcommand{\ecor}{\end{corr}}
\newcommand{\edeff}{\end{deff}}

\newcommand{\eprop}{\end{proposition}}

\renewcommand{\Pi}{\varPi}

\renewcommand{\epsilon}{\varepsilon}

\begin{document}

\title%[Global existence for wave equations]
{A natural lower bound for the size of nodal sets}
\thanks{The authors were supported in part by  NSF grants DMS-0969745 and DMS-1069175.}

\author{Hamid Hezari}
\author{Christopher D. Sogge}
\address{Department of Mathematics, M.I.T., Cambridge, MA 02139}
\address{Department of Mathematics,  Johns Hopkins University,
Baltimore, MD 21218}

\maketitle

%\newsection{Introduction}

The purpose of this brief note is to prove a natural lower bound for the $(n-1)$-dimensional Hausdorff measure of nodal sets of eigenfunctions.  To wit:

\begin{theorem}\label{maintheorem}  Let $(M,g)$ be a compact manifold of dimension $n$ and
$e_\lambda$ an eigenfunction satisfying 
$$-\Delta_g e_\lambda = \lambda e_\lambda, \, \, \text{and } \, \, 
\int_{M}|e_\lambda|^2 \, dV_g=1.$$
Then if $Z_\lambda=\{x\in M: \, e_\lambda(x)=0\}$ is the nodal set and $|Z_\lambda|$ its
$(n-1)$-dimensional Hausdorff measure, we have
\begin{equation}\label{main}
\lambda^{\frac12}\left(\, \int_M |e_\lambda|\, dV_g\, \right)^2 \le C |Z_\lambda|, \quad \lambda \ge1,
\end{equation}
for some uniform constant $C$.  Consequently,
\begin{equation}\label{5}
\lambda^{\frac{3-n}4}\lesssim |Z_\lambda|, \quad \lambda\ge1.
\end{equation}
\end{theorem}

Inequality \eqref{5} follows from \eqref{main} and the lower bounds in \cite{SZ}
\begin{equation}\label{3}
\lambda^{\frac{1-n}8} \lesssim \int_M |e_\lambda| \, dV_g.
\end{equation}
The lower bound \eqref{5} is due to Colding and Minicozzi~\cite{CM}.  Yau~\cite{Y} conjectured
that $\lambda^{\frac12}\approx |Z_\lambda|$.  This lower bound $\lambda^{\frac12}
\lesssim |Z_\lambda|$ was verified in the 2-dimensional case by Br\"uning~\cite{B} and independently by Yau (unpublished).  The bounds in \eqref{5} 
seem to be the best known ones for higher dimensions, although Donnelly and Fefferman~\cite{DF1}-\cite{DF2} showed that, as conjectured, $|Z_\lambda|\approx \lambda^{\frac12}$, if $(M,g)$ is assumed to be real analytic.  

The first ``polynomial type'' lower bounds appear to be due to to Colding and Minicozzi~\cite{CM} and Zelditch and the second author~\cite{SZ} (see also \cite{M}).  As we shall point out inequality~\eqref{main} cannot be improved and it to some extent unifies the approaches in \cite{CM} and
\cite{SZ}.  As was shown in \cite{SZ}, the $L^1$-lower bounds in \eqref{3} follow from
H\"older's inequality and the 
$L^p$ eigenfunction estimates of the second author~\cite{S2}
for the range where $2<p\le \frac{2(n+1)}{n-1}$.  These too cannot be improved,
but it is thought better $L^p$-bounds hold for a typical eigenfunction or if one makes geometric
assumptions such as negative curvature (cf. \cite{SZ2}-\cite{SZ3}).  
Thus, it is natural to expect to be able to improve \eqref{3} and hence the lower bounds
\eqref{5} for all eigenfunctions on manifolds with negative curvature, or for ``typical" eigenfunctions on any manifold.  Of course, Yau's conjecture that $|Z_\lambda|\approx \lambda^{\frac12}$ would be the ultimate goal, but understanding when \eqref{3} 
can be improved is a related problem of independent interest.

\bigskip

Let us now turn to the proof of Theorem~\ref{maintheorem}.  We shall use an identity
from the recent work of the second author and Zelditch~\cite{SZ}:
\begin{equation}\label{1}
\int_M |e_\lambda| \, (\Delta_g+\lambda)f \, dV_g = 2 \int_{Z_\lambda} |\nabla_g e_\lambda| \, f
\, dS_g,
\end{equation}
Here $dS_g$ is the Riemannian surface measure on $Z_\lambda$,
and  
$\nabla_g$ is the gradient coming from the metric and $|\nabla_g u|$ is the norm
coming from the metric, meaning that in local coordinates
\begin{equation}\label{quad}|\nabla_g u|^2_g = \sum_{jk=1}^n g_{jk}(x)\partial_j u \partial_k u.
\end{equation}
Identity ~\eqref{1} follows from the
Gauss-Green formula and a related earlier identity was proved by Dong~\cite{D}.

As in \cite{HW}, if
we take $f\equiv1$ and apply Schwarz's inequality we get
\begin{equation}\label{2}
\lambda \int_M |e_\lambda| \, dV_g \le 2|Z_\lambda|^{1/2} \, \left(\, 
\int_{Z_\lambda}|\nabla_g e_\lambda|^2 \, dS_g\, \right)^{1/2}.
\end{equation}
Thus we would have \eqref{main} if we could prove that the energy of $e_\lambda$ on its nodal set satisfies the natural bounds
\begin{equation}\label{4}
\int_{Z_\lambda}|\nabla_g e_\lambda|^2 \, dS_g \lesssim \lambda^{\frac32}.
\end{equation}

We shall do this by choosing a different auxiliary function $f$.
This time we want to use
\begin{equation}\label{6}
f=\big(\, 1+\lambda e_\lambda^2 + |\nabla_g e_\lambda|^2_g \, \bigr)^{\frac12}.
\end{equation}
If we plug this into \eqref{1} we get that
\begin{equation*}
2\int_{Z_\lambda}|\nabla_g e_\lambda|^2_g dS_g
\le \int_M |e_\lambda| \, (\Delta_g+\lambda)\big(\,1+ \lambda e_\lambda^2 + |\nabla_g e_\lambda|^2 \, \bigr)^{\frac12} \, dV_g.
\end{equation*}
Since we have the $L^2$-Sobolev bounds
\begin{equation}\label{7}\|  e_\lambda\|_{H^s(M)}=O(\lambda^{\frac{s}2}),
\end{equation}
it is clear that
$$\lambda \int_M |e_\lambda| \, \big(\,1+ \lambda e_\lambda^2 + |\nabla_g e_\lambda|^2_g \, \bigr)^{\frac12}\, dV_g = O(\lambda^{\frac32}),
$$
and thus to prove \eqref{4}, it suffices to show that
\begin{equation}\label{8}
\int_M |e_\lambda| \, \Delta_g\big(\, 1+\lambda e_\lambda^2 + |\nabla_g e_\lambda|^2_g \, \bigr)^{\frac12}
\, dV_g = O(\lambda^{\frac32}).
\end{equation}
To prove this we first note that
$$
\partial_k \big(\, 1+\lambda e_\lambda^2 + |\nabla_g e_\lambda|^2_g \, \bigr)^{\frac12}
=\frac{\lambda e_\lambda \partial_k e_\lambda
+ \frac12\partial_k |\nabla_g e_\lambda|_g^2
%\sum_l \partial_k \partial_l e_\lambda \partial_l e_\lambda
}
{\big(\, 1+ \lambda e_\lambda^2 + |\nabla_g e_\lambda|^2 \, \bigr)^{\frac12}},
$$
from this and \eqref{7} we deduce that
$$\int_M |e_\lambda| \left| \, \nabla_g \big(\, 1+\lambda e_\lambda^2 + |\nabla_g e_\lambda|^2 \, \bigr)^{\frac12} \, \right| \, dV_g = O(\lambda).$$
This means that the contribution of the first order terms of the Laplace-Beltrami operator
(written in local coordinates) to 
\eqref{8} are better than required, and so it suffices to show that in a 
compact subset $K$ of a
local coordinate patch we have
\begin{equation}\label{9}
\int_K |e_\lambda| \, 
\Bigl|  \partial_j \partial_k \big(\, 1+\lambda e_\lambda^2 + |\nabla_g e_\lambda|^2 \, \bigr)^{\frac12}
\, \Bigr| \, dV_g = O(\lambda^{\frac32}).
\end{equation}
A calculation shows that
$\partial_j\partial_k \big(\, \lambda e_\lambda^2 + |\nabla_g e_\lambda|^2 \, \bigr)^{\frac12}$ equals
\begin{multline*}
-\frac{\bigl(\lambda e_\lambda \partial_j e_\lambda + 
\frac12\partial_j |\nabla_g e_\lambda|^2_g
%\sum_l \partial_j\partial_l e_\lambda 
%\partial_l e_\lambda
\bigr)
\bigl(\lambda e_\lambda \partial_k e_\lambda + 
\frac12\partial_k |\nabla_g e_\lambda|^2_g
%\sum_l \partial_k\partial_l e_\lambda 
%\partial_l e_\lambda
\bigr)
}
{ \big(\,1+ \lambda e_\lambda^2 + |\nabla_g e_\lambda|^2 \, \bigr)^{\frac32} }
\\
+
\frac{\lambda\partial_j e_\lambda \partial_ke_\lambda+
\lambda e_\lambda \partial_j\partial_k e_\lambda
+
\frac12\partial_j\partial_k |\nabla_g e_\lambda|^2_g
%\sum_l \partial_j \partial_k\partial_l e_\lambda \partial_l e_\lambda
%+\sum_l \partial_k\partial_l e_\lambda \partial_j\partial_l e_\lambda
}
{ \big(\,1+ \lambda e_\lambda^2 + |\nabla_g e_\lambda|^2 \, \bigr)^{\frac12} }
.
\end{multline*}
If $|D^m f|=\sum_{|\alpha|=m}|\partial^\alpha f|$, then by \eqref{quad}
$$\partial_k |\nabla_g e_\lambda|^2 =O(|D^2e_\lambda| \, |De_\lambda| +|De_\lambda|^2),$$
and
$$\partial_j\partial_k |\nabla_g e_\lambda|^2_g 
=O(|D^3 e_\lambda| \, |De_\lambda| + |D^2e_\lambda|^2 +
 |D^2 e_\lambda|\,
|De_\lambda|+|De_\lambda|^2 ).$$
Therefore,
\begin{multline*}
\partial_j\partial_k \big(\, \lambda e_\lambda^2 + |\nabla_g e_\lambda|^2 \, \bigr)^{\frac12}
=O\left(\frac{\lambda^2 |e_\lambda|^2 |D e_\lambda|^2+|D^2 e_\lambda|^2
|D e_\lambda|^2 +|De_\lambda|^4}{
 \big(\,1+ \lambda e_\lambda^2 + |\nabla_g e_\lambda|^2 \, \bigr)^{\frac32}}\right)
 \\
 +O\left( 
 \frac{\lambda |D e_\lambda|^2 + \lambda |e_\lambda| \, |D^2 e_\lambda|
 +|D^3 e_\lambda| \, |D e_\lambda|+|D^2 e_\lambda|^2+|D^2e_\lambda| \, |De_\lambda|
 +|De_\lambda|^2}
 { \big(\, 1+\lambda e_\lambda^2 + |\nabla_g e_\lambda|^2 \, \bigr)^{\frac12} }
 \right).
 \end{multline*}
This implies that the integrand in the left side of \eqref{9} is dominated by
\begin{multline*}
\Bigl( \, 
\lambda^{\frac12}|D e_\lambda|^2+\lambda^{-\frac12}|D^2 e_\lambda|^2
\, +|De_\lambda|^2 \Bigr)
\\
+
\Bigl( \, 
\lambda^{\frac12}|D e_\lambda|^2+ \lambda^{\frac12} |e_\lambda| \, |D^2 e_\lambda|
+|e_\lambda| \, |D^3 e_\lambda| +\lambda^{-\frac12}|D^2 e_\lambda|^2
+|D^2e_\lambda| \, |e_\lambda| + |De_\lambda| \, |e_\lambda|
\, \Bigr),
\end{multline*}
leading to \eqref{9} after applying \eqref{7}.  \qed

\noindent{{\bf Remarks:}  
\begin{itemize}
\item
We could also have taken $f$ to be $(\lambda+\lambda e_\lambda^2
+|\nabla_g e_\lambda|^2)^{\frac12}$ and obtained the same upper bounds, but there does  not seem to be any advantage to doing this.  
\item Inequality \eqref{main} cannot be improved.  There are many cases when the $L^1$ and 
$L^2$-norms of eigenfunctions are comparable.  For instance, for the sphere the zonal functions have this property and it is easy to check that their nodal sets satisfy $|Z_\lambda|\approx \lambda^{\frac12}$, which means that for zonal functions \eqref{main} cannot be improved.
\item There are many cases where inequality \eqref{main} can be improved.  For instance,
the $L^2$-normalized highest weight spherical harmonics $Q_k$ have eigenvalues $\lambda
=\lambda_k \approx k^2$, and $L^1$-norms $\approx k^{-\frac{n-1}4}$ (see e.g., \cite{S1}).  This means that for
the highest weight spherical harmonics the left side is proportional to $\lambda^{\frac{3-n}4}$
even though here too $|Z_\lambda|\approx \lambda^{\frac12}$. 
Similarly, the highest weight spherical harmonics saturate \eqref{4}.  It is because of functions like
the highest weight spherical harmonics that the current techniques only seem to yield
\eqref{5}.  Note that inequality \eqref{5} gives the correct lower bound in the trivial case
where the dimension $n$ is one.  As the dimension increases, the bound gets worse and worse
due to the fact that \eqref{3} is saturated by functions like the highest weight spherical harmonics
(``Gaussian beams'') whose mass is supported on a $\lambda^{-\frac14}$ neighborhood of a geodesic and the volume of such a tube decreases geometrically as $n$ increases.  (See \cite{bourgainef}
and \cite{Sokakeya} for related work on this phenomena.)
\end{itemize}

\medskip
\noindent {\bf Acknowledgments}  The authors wish to thank W. Minicozzi and S. Zelditch for several helpful and interesting discussions.

\end{document}